\documentclass[11pt,a4paper]{article}

\usepackage{amsmath,amsfonts,amssymb}
\usepackage{color}
\usepackage{a4wide}
\usepackage{amscd}
\usepackage{url}
\newtheorem{theorem}{Theorem}
\newtheorem{conjecture}[theorem]{Conjecture}

\def\R{\mathbb R}
\def\Z{\mathbb Z}

\def\Prob{\mathrm{Prob}}
\def\div{\mathop{\mathrm{div}}}
\DeclareMathOperator\grad{grad}
\def\X{\mathcal X}

\def\M{\mathcal M}
\DeclareMathOperator\Ent{Ent}
\DeclareMathOperator\Entmix{Ent_{mix}}
\def\Lebesgue{\mathcal L}
\def\E{\mathcal E} 
\def\F{\mathcal F}
\DeclareMathOperator*\argmin{arg\,min}
\let\ds\displaystyle
\def\un#1{\,\mathrm{#1}}
\def\Y{\mathcal Y}
\def\H{\mathcal H}
\def\P{\mathcal P}

\begin{document}

\title{Large deviations and gradient flows}

\author{Stefan Adams\footnote{Mathematics Institute, University of Warwick}\and Nicolas Dirr\footnote{Cardiff School of Mathematics, Cardiff University}\and Mark Peletier\footnote{Department of Mathematics and Computer Sciences and Institute for Complex Molecular Systems, Technische Universiteit Eindhoven} \and Johannes Zimmer\footnote{Department of Mathematical Sciences, University of Bath}}
\maketitle

\begin{abstract}
In recent work~\cite{AdamsDirrPeletierZimmer11} we uncovered intriguing connections between Otto's   characterisation of diffusion as entropic gradient   flow~\cite{Otto01} on one hand and large-deviation principles   describing the microscopic picture (Brownian motion) on the   other. In this paper, we sketch this connection, show how it generalises to a   wider class of systems, and comment on consequences and implications. 
  
Specifically, we connect macroscopic gradient flows with large deviation principles, and point out the potential of a bigger picture emerging: we indicate that in some non-equilibrium situations, entropies and thermodynamic free energies can be derived via large deviation principles. The approach advocated here is different from the established hydrodynamic limit passage but extends a link that is well known in the equilibrium situation.
\end{abstract}

\section{Introduction}

For systems in equilibrium, it is well known that the roles of energy and entropy can be understood rigorously in terms of large-deviation principles. We describe two examples below. Recently, we showed how large-deviation principles also allow us to understand the role of entropy in a specific \emph{non-equilibrium} system~\cite{AdamsDirrPeletierZimmer11}: the large-deviation behaviour of a system of independent Brownian particles connects rigorously to the entropy gradient-flow structure of the diffusion equation. We explain this connection in Section~\ref{subsec:just-diff}. 

The aim of this paper is to take this connection two steps further. The first step is to extend the  connection of~\cite{AdamsDirrPeletierZimmer11}, which was studied in a discrete-time context, to the case of continuous time. The second step is to discuss a variety of examples that illustrates the breadth of this phenomenon, and suggest a general principle that might hold across a wide range of systems. 

\medskip
 
In equilibrium systems, the connection is as follows. 
Let $X_i$ ($i=1,2,\dots$) be independent and identically distributed stochastic variables with distribution $\mu$ on a state space $\X$. We think of the $X_i$ as positions of particles in the space $\X$, so that their concentration is given by the \emph{empirical measure} $\rho_n := \frac1n  \sum_{i=1}^n \delta_{X_i}$. Sanov's theorem (e.g.,~\cite[Sec.~6.2]{DemboZeitouni98}) states that the random measure~$\rho_n$ satisfies the \emph{large-deviation principle}
\begin{equation}
\label{ldp:entropy}
\mathrm{Prob}(\rho_n \approx \rho)\sim \exp [-n I(\rho)],
\qquad\text{as } n\to\infty,
\end{equation}
where the \emph{rate function} $I\geq0$ is the \emph{relative entropy} of $\rho$ with respect to $\mu$, which is
\begin{equation*}
I(\rho) = H(\rho|\mu) := 
\begin{cases}
\int f\log f\, d\mu & \text{if }\rho\ll \mu \text{ and } \rho = f\mu,\\
+\infty & \text{otherwise},
\end{cases}
\end{equation*}
This property illustrates how the relative entropy $H(\rho|\mu)$ characterises the probability of observing a state $\rho$: higher relative entropy means smaller probability, as described by~\eqref{ldp:entropy}. It also provides a rigorous version of the well-known thermodynamic principle that a system aims to maximise its entropy (which corresponds to minimising $H(\rho|\mu)$, since the physical entropy carries the opposite sign). For in the limit of large $n$, the characterisation~\eqref{ldp:entropy} gives vanishing probability to all states $\rho$ except those for which $I(\rho) = 0$; in other words, only the minimisers of $I$ have non-vanishing probability. 
%The proof of this theorem shows that the relative entropy arises from the indistinguishibility of 

This connection between entropy and large-deviation principles extends to systems involving energy. In the appendix we show, for instance, how coupling a system with energy $E$ to a heat bath with temperature $\theta$ changes the rate functional $I$ to the \emph{free energy} $\mathcal F(\rho) := H(\rho|\mu) + (k\theta)^{-1} E(\rho) + \text{constant}$:
\begin{equation}
\label{ldp:free-energy}
\mathrm{Prob}(\rho_n \approx \rho)\sim \exp [-n \mathcal F(\rho)],
\qquad\text{as } n\to\infty.
\end{equation}
In the same way as~\eqref{ldp:entropy} explains why relative entropy is minimized, \eqref{ldp:free-energy} explains why systems coupled to a heat bath minimize their free energy: when $n$ is large, only states $\rho$ with near-minimal free energy $\F(\rho)$ will have finite probability.

\medskip

As mentioned above, the central aim of this paper is to show how this connection between entropy and free energies on one hand and large-deviation principles on the other extends into the realm of non-equilibrium systems. We restrict our focus to the important class of \emph{gradient flows}, where this connection explains many aspects of these systems. Since the entropy appears as the driving force of the process, we will occasionally call this functional ``energy'' to conform with the standard terminology for gradient flows.

The general philosophy is illustrated by the diagram below. 
\makeatletter
\atdef@ O#1O#2O{\CD@check{O..O..O}{\llap{$\m@th\vcenter{\hbox 
 {$\scriptstyle#1$}}$}\Big\updownarrow 
 \rlap{$\m@th\vcenter{\hbox{$\scriptstyle#2$}}$}&&}} 
\atdef@+#1+#2+{\ampersand@
  \ifCD@ \global\bigaw@\minCDarrowwidth \else \global\bigaw@\minaw@ \fi
  \setboxz@h{$\m@th\scriptstyle\;\;{#1}\;$}%
  \ifdim\wdz@>\bigaw@ \global\bigaw@\wdz@ \fi
  \@ifnotempty{#2}{\setbox\@ne\hbox{$\m@th\scriptstyle\;\;{#2}\;$}%
    \ifdim\wd\@ne>\bigaw@ \global\bigaw@\wd\@ne \fi}%
  \ifCD@\enskip\fi
    \mathrel{\mathop{\hbox to\bigaw@{\leftrightarrowfill@\displaystyle}}%
      \limits^{#1}\@ifnotempty{#2}{_{#2}}}%
  \ifCD@\enskip\fi \ampersand@}
\makeatother
\def\ntext#1{\text{\normalsize #1}}
\begin{equation}
\label{diag:commuting}
\begin{CD} 
\substack{\textstyle\text{\emph{dynamic} rate functional}\\[1\jot]\textstyle I \text{ or }I_h}@+\text{\quad this paper\quad}++
 \substack{\textstyle\text{gradient-flow structure}\\[\jot]\textstyle\ntext{$J$ or $I_h$}}\\ 
@O\substack{\text{large-deviation principle}\\n\to\infty}OO 
  @OOO\\ 
\text{stochastic $n$-particle system} @>\text{continuum limit}>n\to\infty> \text{continuum evolution equation} 
\end{CD}
\end{equation}
The bottom row in this diagram is the classical connection between a stochastic $n$-particle system and its hydrodynamic limit: the typical case is that as $n\to\infty$, the particle system becomes deterministic, and the empirical measure of the particle system converges to the solution of the (deterministic) continuum equation. Note that this statement concerns only the typical behaviour of the particle system; large deviations are not captured.

In the left-hand column, a large-deviation principle characterises the behaviour in the limit $n\to\infty$ in a different manner, in terms of a functional $I$ or $I_h$ of the \emph{time-dependent} system, as we shall see below.
The right-hand column is the connection between an evolution equation and the corresponding gradient-flow structure, when it exists. 

The central statement of this paper is the double-headed arrow at the top. It provides a connection between representations with more information on both sides: on the left-hand side, the rate functional contains more information than just the most probable behaviour, and on the right-hand side, the gradient-flow structure is an additional structure on top of the equation itself. 

In the following sections, we illustrate the double-headed arrow in a number of concrete examples, first in the discrete-time approximation (Section~\ref{sec:discrete}) and then in continuous time (Section~\ref{sec:continuous}). Section~\ref{sec:Other-gradient-flows} generalises the argument to non-quadratic dissipations. Since the implications of this connection are best appreciated once one has an overview of the breadth of the phenomenon, we postpone most of the discussion of the consequences to Section~\ref{sec:discussion}. 
\medskip

The mathematical results described in this paper are not new, and mostly due to other authors, such as Freidlin \& Wentzell~\cite{FreidlinWentzell98}, Dawson \& G\"artner~\cite{DawsonGartner87,DawsonGartner89}, Feng \& Kurtz~\cite{FengKurtz06}, Kipnis, Olla, \& Varadhan~\cite{KipnisOllaVaradhan89} and others. Instead, we see the novelty of this paper in extracting from these results the suggestion of a general principle connecting the broad class of gradient flows with large deviations of stochastic processes. A particularly interesting aspect of this connection is that thermodynamic quantities are derived in a  non-equilibrium context.

\section{The Wasserstein metric}
\label{sec:Wasserstein}

Much of this paper centres on the Wasserstein metric and Wasserstein gradient flows. 
The (quadratic) \emph{Wasserstein distance} between two probability measures $\rho_0$ and $\rho_1$ with finite second moments is~\cite{Villani03}
\begin{equation}
\label{def:d}
d(\rho_0,\rho_1)^2 = \inf_q \int_{\R^d\times\R^d} |x-y|^2 \, q(dxdy),
\end{equation}
where the infimum is taken over all $q$ with marginals $\rho_0$ and $\rho_1$, i.e., over all $q$ satisfying 
\begin{equation*}
\text{for any $A\subset \R^d$}, \qquad
q(A\times \R^d) = \rho_0(A) \qquad\text{and}\qquad
q(\R^d\times A) = \rho_1(A).
\end{equation*} 

We  also need an incremental version of the Wasserstein distance. The Brenier-Benamou formula~\cite{BenamouBrenier00} gives
 an alternative formulation of $d$ as an infimum of curves of measures $t\mapsto \rho(t)$ such that $\rho(0) = \rho_0$ and $\rho(1)=\rho_1$:
\begin{equation}
\label{def:W-distance}
d(\rho_0,\rho_1)^2 = \inf_{\rho\colon[0,1]\to\mathcal{M}_1(\mathbb{R}^d)} \int_0^1 \|\partial_t \rho(t)\|^2_{\rho(t),*}\, dt.
\end{equation}
Here the local norm $\|\cdot\|_{\rho,*}$ at a given point $\rho$ is derived from an inner product (a local metric tensor) formally given by
\begin{equation}
\label{def:g}
(s_1,s_2)_{\rho,*} := \int_{\R^d} \rho(x) \nabla p_1(x)\cdot\nabla p_2(x)\, dx,
\end{equation}
where   $\nabla$ is the usual gradient in  $\R^d$, and the $p_i$ solve the equation $\div (\rho\nabla p_i) = s_i$ in $\R^d$ (see~\cite{DawsonGartner87,KipnisOlla90} or~\cite[Sec.~9.4]{FengKurtz06} for a rigorous definition). 
\medskip

A \emph{Wasserstein gradient flow} is a gradient flow of an energy $\E$ with respect to the Wasserstein metric structure. A curve of measures $t\mapsto \rho(t)$ is a solution of such a gradient-flow equation if its time derivative $\partial_t\rho$, in the sense of distributions, satisfies
\begin{equation}
\label{def:innerproductstar}
(\partial_t\rho(t),s_2)_{\rho(t),*} = -\int_{\R^d} \frac{\delta\E}{\delta \rho}(\rho(t))\,s_2\, dx
\qquad \text{for all }s_2\text{ and all }t>0,
\end{equation}
where $\delta \E/\delta\rho$ is the variational derivative of~$\E$.
A straightforward calculation shows that this is equivalent to the equation
\begin{equation}
\label{eq:WGF}
\partial_t \rho = \div \rho\nabla \Bigl(\frac{\delta\E}{\delta \rho}\Bigr).
\end{equation}
By analogy with gradients in Riemannian geometry, this suggests to define the \emph{Wasserstein gradient} of a functional $\E$ as
\begin{equation}
\label{formula:GF-Wass}
\grad_W \E(\rho) := -\div \rho\nabla\Bigl( \frac{\delta \E}{\delta \rho}\Bigr).
\end{equation}

Below we shall also use more general versions of this structure. Replacing $\rho$ above by a general diffusion matrix $D(\rho)$, we define
\begin{equation}
\label{def:gD}
(s_1,s_2)_{D(\rho),*} := \int_{\R^d} D(\rho(x)) \nabla p_1(x)\cdot\nabla p_2(x)\, dx,
\qquad\text{where} \qquad 
s_i = \div D(\rho)\nabla p_i.
\end{equation}
Repeating the construction above, it follows that the \emph{$D$-Wasserstein gradient} of a functional~$\E$ is characterised by the equation
\begin{equation}
\label{formula:GF}
\partial_t \rho = \div D(\rho)\nabla\Bigl( \frac{\delta \E}{\delta \rho}\Bigr).
\end{equation}

%With this definition the diffusion equation~\eqref{eq:diffusion} can indeed be written as the Wasserstein gradient-flow of $\Ent$ (with $D(\rho)=\rho$), i.e.,
%\begin{equation*}
%\partial_t \rho = -\grad_{W,\rho} \Ent(\rho).
%\end{equation*}

\medskip

Gradient flows have natural time-discrete approximations, constructed in an iterative manner: 
\begin{align}
&\textit{For given approximation $\rho_{k-1}$ at time $(k-1)h$, choose $\rho_k$ at time $kh$}\notag\\
&\qquad \textit{as minimiser of the functional }\quad \rho\mapsto  \frac1{2h}d(\rho,\rho^{k-1})^2 + \E(\rho ).
\label{approx:JKO}
\end{align}
This is essentially a backward-Euler discretisation, as can be recognised by comparing it with the $\R^d$-gradient-flow $\dot x = -\nabla E(x)$. For this equation the backward-Euler discretisation is constructed by solving
\begin{equation*}
 \frac1h (x_k-x_{k-1}) = -\nabla E(x_k),
\end{equation*}
for $x_k$, which is equivalent to minimising 
\begin{equation}
\label{GF:l2a}
x\mapsto \frac1{2h} |x-x_{k-1}|^2 + E(x).
\end{equation}
Note the similarity between~\eqref{GF:l2a} and~\eqref{approx:JKO}: in both expressions the first term measures the distance between old and new states, while the second term favours a reduction of the functional~$\E$ respectively $E$.

\section{Discrete time}
\label{sec:discrete}

We can now formulate the first example. 

\subsection{A system of independent Brownian particles}
\label{subsec:just-diff}

We consider $n$ independent Brownian particles $X_{n,i}(t)$ in $\R^d$, with deterministic initial positions $X_{n,i}(0) = x_{n,i}$, each hopping to a new position $X_{n,i}(h)$ at time $h>0$ with a Gaussian probability with mean $x_{n,i}$ and variance\footnote{In this paper, we consider Brownian particles with generator $\Delta$, rather than $(1/2)\Delta$, and therefore the transition kernel is $(4\pi h)^{-d/2}\exp{-|x-y|^2/4h}$.}~$2h$. 

As in the equilibrium case discussed above, we describe this system by the empirical measure $\rho_n(t):= \tfrac1n \sum_{i=1}^n \delta_{X_{n,i}(t)}$ at a given time $t$, and we assume that the initial measure $\rho_n(0)$ converges to a given measure $\rho^0$ as $n\to\infty$. 
In the limit of large~$n$, the probability of this jump process attaining any $\rho^1$ at time $t=h$ is again characterised in terms of a {large-deviation principle},
\begin{equation}
\label{ldp:diffusioneq}
\mathrm{Prob}(\rho_n(h) \approx \rho^1)\approx \exp [-n I_h(\rho^1)],
\end{equation}
where the {rate functional} $I_h$ has an explicit expression that can be derived from Stirling's formula (see~\cite{AdamsDirrPeletierZimmer11} for the expression; in~\cite{AdamsDirrPeletierZimmer11}, $I_h$ is  only the limit of a sequence of rate functionals, but can be shown to be a rate functional in its own right~\cite{Leonard07TR,PeletierRenger11TR}). 

The main result of~\cite{AdamsDirrPeletierZimmer11} is that 
\begin{equation}
\label{eq:ADPZ}
  I_h \approx K_h \qquad \text{as }h\to0,
\end{equation}
where 
\begin{equation}
\label{def:K_h}
K_h(\rho^1;\rho^0) := \frac1{4h} d(\rho^0,\rho^1)^2 + \frac12 \Ent(\rho^1) - \frac12 \Ent(\rho^0).
\end{equation}
Here $d$ is the Wasserstein distance defined above, and 
\begin{equation*}
\Ent(\rho):= H(\rho|\Lebesgue) = \begin{cases}
\int f\log f\, dx & \text{if }\rho\ll\Lebesgue \text{ and }\rho = f\Lebesgue,\\
+\infty & \text{otherwise},
\end{cases}
\end{equation*}
is the {relative entropy} of $\rho$ with respect to the Lebesgue measure $\Lebesgue$. The rigorous formulation of~\eqref{eq:ADPZ} is a Gamma-convergence result of $I_h$ to $K_h$ after both have been desingularised. 

\medskip

The functional $K_h$ has the same form as the functional in~\eqref{approx:JKO}, since the term $\Ent(\rho^0)/2$ does not influence the minimisation with respect to $\rho^1$. Therefore the time-discrete approximation that one constructs with this $K_h$ is an approximation of the Wasserstein gradient flow of the entropy $\Ent$, which is the  diffusion equation~\cite{JordanKinderlehrerOtto98}
\begin{equation}
\label{eq:diffusion}
\partial_t\rho = \Delta \rho 
\qquad \text{in }\R^d.
\end{equation}

This is the connection referred to above: the large-deviation behaviour of the system of particles is represented by the rate functional $I_h$, and this functional is asymptotically equal to the functional $K_h$ that defines the gradient-flow formulation of the diffusion equation. The approximation result~\eqref{eq:ADPZ} therefore creates a link between the gradient-flow structure of the deterministic limit equation on one hand and the large-deviation behaviour of the system of particles on the other. The same result can be shown for Gaussian measures on the real line~\cite{Dirr:11a}. In the rest of this paper we shall see many more versions of such connections. 

\paragraph{Consequences}
While most of the discussion is deferred to Section~\ref{sec:discussion}, we mention here a few consequences of the fact~\eqref{eq:ADPZ} that  the large-deviation rate functional $I_h$ and the constructing functional $K_h$ of the gradient flow are equal in the limit $h\to0$.

First, the construction of a time-discrete approximation~\eqref{approx:JKO} to the diffusion equation~\eqref{eq:diffusion} was motivated in~\cite{JordanKinderlehrerOtto98} by analogy with the backward-Euler discretisation~\eqref{GF:l2a}. This is an indirect and purely mathematical motivation, which explains neither the reason for the appearance of the entropy and the Wasserstein distance in $K_h$, nor the reason for minimising just this combination. 

The connection between $K_h$ and~$I_h$, however, gives a direct motivation. By~\eqref{ldp:diffusioneq}--\eqref{eq:ADPZ},  $K_h(\rho;\rho^0)$ is a measure of the likelihood of observing a state $\rho$ after time $h$. For large $n$, the characterisation~\eqref{ldp:diffusioneq} implies that only the global minimiser of~$I_h$, and therefore of $K_h$, is observed with non-vanishing probability. The stochastic minimisation~\eqref{ldp:diffusioneq} of $I_h$ thus becomes converted into an absolute minimisation of $K_h$.

\medskip
Secondly, in the limit $h\to0$, the \emph{proof} that $I_h\approx K_h$ explains the origin of the two terms of $K_h$. The entropy arises from the indistinguishibility of the particles after transforming to an empirical measure. The origin of the Wasserstein cost functional $|x-y|^2$ in~\eqref{def:d} can be traced back to the exponent of the term $e^{-|x-y|^2/4h}$ in the Gaussian transition probability of the Brownian particles. We return to this issue in Section~\ref{sec:discussion}.

\section{Continuous time}
\label{sec:continuous}

The construction in the previous section is discrete in time: the rate function $I_h$ describes the probability distribution of the state $\rho_n(h)$ at time $h>0$. A continuous-time large-deviation principle, where one considers deviations from a whole path of empirical measures for a fixed terminal time,  provides a different kind of insight, and may be even closer to the gradient-flow formulation. We start with some preliminaries.

\subsection{An alternative formulation of the gradient-flow structure}
\label{subsec:alt-GF-formulation}

In a formal sense, Wasserstein gradient flows and many others can  be written in the form
\begin{equation}
\label{def:M-GF}
\partial_t \rho =- M_\rho\frac{\delta\E}{\delta \rho},
\end{equation}
where $\E$ is the `energy' functional driving the evolution, and $M_\rho$ a $\rho$-dependent symmetric mapping\footnote{This way of writing the gradient flow highlights the fact that a gradient flow is an instance of a GENERIC evolution, in which the conservative evolution term is absent~\cite{Oettinger05}.}. In the case of Wasserstein gradient flows, for instance, 
\begin{equation*}
M_\rho\xi = -\div \rho\nabla \xi,
\end{equation*}
as follows by comparing~\eqref{eq:WGF} with~\eqref{def:M-GF}.
Taking this case of Wasserstein gradient flow as an example, 
we shall encounter  the equation~\eqref{def:M-GF} in a different form, connected to the functional $J$ given by
\begin{equation}
\label{def:J}
J(\rho) := \E(\rho(T)) - \E(\rho(0)) + \frac12 \int_0^T\left[\|\partial_t \rho\|_{\rho,*}^2 + \Bigl\|-\frac{\delta \E}{\delta\rho}\Bigr\|^2_{\rho}\right]\, dt,
\end{equation}
where 
\begin{equation*}
\|\xi\|_{\rho}^2  := \int_{\R^d} \xi \, M_\rho\,\xi\, dx = \int_{\R^d} \rho\,|\nabla\xi|^2
\end{equation*}
and the norm $\|\cdot\|^2_{\rho,*}$ is the norm defined in~\eqref{def:g}. The norms $\|\cdot\|_\rho$ and $\|\cdot\|_{\rho,*}$ are dual norms, and $\|\cdot\|_{\rho,*}$ has the alternative characterisation
\begin{equation*}
\|s\|_{\rho,*} := \sup_{\xi\not=0} \frac{\ds\int_{\R^d} s\xi\, dx}{\|\xi\|_{\rho}}.
\end{equation*}
By writing the energy difference $\E(\rho(T))-\E(\rho(0))$ as 
\begin{equation*}
\E(\rho(T))-\E(\rho(0)) = \int_0^T \int_{\R^d} \frac{\delta \E}{\delta \rho}\, \partial_t\rho \, dx dt 
= \int_0^T \Bigl(M_\rho\frac{\delta \E}{\delta \rho}, \partial_t \rho\Bigr)_{\rho,*}\, dt , 
\end{equation*}
using the inner product defined in~\eqref{def:innerproductstar},
the functional~$J$ in~\eqref{def:J} can now be written as 
\begin{equation*}
J(\rho) \;=\; \frac12\int_0^T \left\|\partial_t \rho + M_\rho\frac{\delta\E}{\delta\rho}\right\|_{\rho,*}^2\, dt.
\end{equation*}
This expression shows that $J$ is non-negative. It also implies that if $\rho$ satisfies $J(\rho)=0$, then  equation~\eqref{def:M-GF} holds at almost each time $0< t< T$; therefore 
\begin{equation}
\label{equiv:connection}
\text{$\rho$ is a Wasserstein gradient flow of $\E$} \quad\Longleftrightarrow \quad J(\rho)=0.
\end{equation}
In the examples of this paper, $J$ is a large-deviation rate functional, and this equivalence  is the connection between the large-deviation behaviour, given by $J$, and the gradient-flow structure of the limiting equation.

\medskip

If we take for the operator $M_\rho$ in~\eqref{def:M-GF} not the Wasserstein operator but a general operator, then we find a similar statement: 
\begin{equation}
\label{equiv:connection-M}
\text{$\rho$ is a solution of the $(\E,M_\rho)$-gradient-flow~\eqref{def:M-GF}} \quad\Longleftrightarrow \quad J_M(\rho)=0,
\end{equation}
where
\begin{equation}
\label{def:JM}
J_M(\rho) := \E(\rho(T)) - \E(\rho(0)) + \frac12 \int_0^T\left[\|\partial_t \rho\|_{M^{-1}_\rho}^2 + \Bigl\|-\frac{\delta \E}{\delta\rho}\Bigr\|^2_{M_\rho}\right]\, dt,
\end{equation}
and the two norms are defined, at least formally, by
\begin{align*}
\|\xi\|_{M_\rho}^2  &:= \int_{\R^d} \xi \, M_\rho\,\xi\, dx,\\
\|s\|_{M^{-1}_\rho} &:= \sup_{\xi\not=0} \frac{\ds\int_{\R^d} s\xi\, dx}{\|\xi\|_{M_\rho}} = \int_{\R^d} s\, M^{-1}_\rho\,s\, dx = \|M^{-1}_\rho\,s\|_{M_\rho}^2.
\end{align*}
We now discuss a number of examples.

%
%In all the examples of this paper $J$ is a large-deviation rate functional, and as we discussed above the curve $\rho$ such that  $J(\rho)=0$ is the only curve that will be observed with non-vanishing probability. This is the general principle of the connection between large-deviation principles and gradient flows, in the case of continuous-time large deviations. 

%The construction in the previous section is discrete in time: the rate function $I_h$ describes the probability distribution of the state $\rho_n(h)$ at time $h>0$. A continuous-time large-deviation principle may provide a different kind of insight, and may be even closer to the gradient-flow formulation. We start with some preliminaries.

\subsection{Continuous-time large deviations for the diffusion equation}

Taking the same system of particles as in Section~\ref{subsec:just-diff}, the continuous-time large-deviation principle for that system of Brownian particles is as follows. Fix a terminal  time $T>0$ and consider the whole path $[0,T]\to\mathcal{M}_1(\mathbb{R}^d) $ of empirical measures $[0,T]\ni t\mapsto\rho_n(t)$. Then the probability that the  \emph{entire curve} $\rho_n(\cdot)$ is close to some other $\rho(\cdot)$ is characterised as~\cite{DawsonGartner87,KipnisOlla90} as a \emph{pathwise large-deviation principle},
\begin{equation*}
\mathrm{Prob}(\rho_n\approx \rho) \sim \exp[-nI(\rho)], 
\end{equation*}
where now
\begin{equation}
\label{I:KipnisOlla}
I(\rho) := \frac12 \int_0^T \left\|\partial_t \rho - \Delta \rho 
\right\|^2_{\rho(t),*}\, dt.
\end{equation}
%Note that $I$ is nonnegative and its minimum of zero is achieved if and only if $\rho$ is a solution of the diffusion equation~\eqref{eq:diffusion}.

This rate function $I$ 
has the structure of $J$ in~\eqref{def:J}.
Using the fact that 
\begin{equation*}
\Delta \rho = \div \rho\nabla \Bigl(\frac{\delta\Ent}{\delta\rho}\Bigr),
\end{equation*}
we find that
\begin{equation*}
I(\rho) = \Ent(\rho(T))-\Ent(\rho(0)) 
 + \frac12\int_0^T \Bigl[ \|\partial_t \rho\|_{\rho,*}^2 + \Bigl\|-\frac{\delta\Ent}{\delta\rho}\Bigr\|_\rho^2\,\Bigr]\, dt.
\end{equation*}
Therefore the Entropy-Wasserstein gradient flow is connected to the large-deviation behaviour of a system of stochastic particles, in the sense of~\eqref{equiv:connection}.
We discuss this further in Section~\ref{sec:discussion}.

\subsection{Diffusive particles with interactions}
\label{subsec:interactions}

We extend the previous example by including interaction of the particles with a background potential $\Psi$ and with each other via an interaction potential~$\Phi$, and modelled by It\^o stochastic differential equations. Specifically, 
we take the microscopic system of $n$ particles to be described by 
\begin{equation}
\label{eq:SDE}
d X_i(t) = -\nabla \Psi(X_i(t))\,dt - \frac1n \sum_{j=1}^n \nabla\Phi(X_i(t)-X_j(t))\, dt + \sqrt 2\, \, dW_i(t),
\end{equation}
where for each $i$, $W_i$ is a Brownian motion in $\R^d$. The hydrodynamic limit of this system is the equation
\begin{equation}
\label{eq:interaction}
\partial_t\rho =  \Delta \rho 
+ \div \rho  \nabla \bigl[\Psi + \rho\ast \Phi\bigr].
\end{equation}

The large-deviation rate functional describing fluctuations of the system is  given by (see~\cite[Theorem~13.37]{FengKurtz06}, and also~\cite{DawsonGartner87} for weakly interacting diffusive particle systems)
\begin{equation}
\label{I:FK}
I(\rho) := \frac12 \int_0^T \Bigl\|\partial_t\rho - \Delta \rho 
- \div \rho \nabla \bigl[\Psi + \rho\ast \Phi\bigr]\Bigr\|^2_{\rho,*}\, dt,
\end{equation}
which again can be written as 
\begin{equation*}
\label{eq:FK-FDT}
I(\rho) 
= \F(\rho(T))- \F(\rho(0)) + \frac12 \int_0^T 
\left[ \|\partial_t \rho\|^2_{\rho,*} 
 + \Bigl\|\frac{\delta\F}{\delta\rho}\Bigr\|^2_\rho
\right] \, dt,
\end{equation*}
where the free energy $\F$ is given by the sum of entropy and potential energy,
\begin{equation}
\label{def:free-energy-F}
\F(\rho) := \Ent(\rho) + \int_{\R^d}\left[\rho\Psi + \frac12 \rho(\rho\ast \Phi)\right].
\end{equation}
Indeed equation~\eqref{eq:interaction} is the Wasserstein gradient flow of the functional~$\F$.

\subsection{The Symmetric Simple Exclusion Process}

The diffusion equation~\eqref{eq:diffusion} is the continuum limit for various stochastic processes, one of which is the system of Brownian particles described above. Here we briefly describe the symmetric simple exclusion process, which has the same limiting equation in a parabolic scaling. However, it has a different large-deviation behaviour, which gives rise to a different gradient flow. 

Consider a periodic lattice $T_n = \{0,1/n, 2/n, \dots (n-1)/n\}$ and its continuum limit, the flat torus $T = \R/\Z$. Each lattice site contains zero or one particle; each particle  attempts to jump from to a neighbouring site with rate $n^2/2$, and they succeed if the target site is empty. We define the configuration $\rho_n\colon T_n\to\{0,1\}$ such that $\rho_n(k/n) = 1$ if there is a particle at site $k/n$, and zero otherwise. For this system the large deviations are characterised by the rate function~\cite{KipnisOllaVaradhan89}
\begin{equation}
\label{I:SSEP}
I(\rho) := \frac12 \int_0^T \|\partial_t \rho - \partial_{xx}\rho\|^2_{\rho(1-\rho),*}\, dt,
\end{equation}
where the norm $\|\cdot\|_{\rho(1-\rho),*}$ is given by~\eqref{def:gD} with  $D(\rho)=\rho(1-\rho)$. This functional can be written as
\begin{equation*}
I(\rho) = \Entmix(\rho(T)) - \Entmix(\rho(0)) 
+ \frac12 \int_0^T \Bigl[ \|\partial_t \rho\|^2_{\rho(1-\rho),*} 
+ \Bigl\|-\frac{\delta\Entmix}{\delta\rho}\Bigr\|^2_{\rho(1-\rho)}\Bigr]\, dt,
\end{equation*}
where the mixing entropy $\Entmix$ is defined as 
\begin{equation*}
\Entmix(\rho) := \int_{\R^d} \bigl[\rho\log\rho + (1-\rho)\log(1-\rho)\bigr].
\end{equation*}
This is true since  $-\partial_{xx} \rho$ is the `$\rho(1-\rho)$'-Wasserstein gradient of $\Entmix$, by
\begin{equation*}
-\partial_{xx}\rho = -\partial_x\left(\rho(1-\rho)\partial_x \log\frac\rho{1-\rho}\right)
= -\partial_x\left(\rho(1-\rho)\partial_x \frac{\delta \Entmix}{\delta\rho} (\rho)\right)
\end{equation*}
(compare this to~\eqref{formula:GF}). Therefore $I$ is of the form~\eqref{def:JM}, with operator
\begin{equation*}
M_\rho\xi := \div \rho(1-\rho) \nabla \xi,
\end{equation*}
and the equation $\partial_t \rho= \partial_{xx}\rho$ is (also) the gradient flow of $\Entmix$ with respect to this `$\rho(1-\rho)$'-Wasserstein structure $\|\cdot\|_{\rho(1-\rho),*}$.

\section{Further generalisations}
\label{sec:Other-gradient-flows}

The arguments of the integrals in~\eqref{def:W-distance}, \eqref{I:KipnisOlla}, \eqref{I:FK}, and~\eqref{I:SSEP} are quadratic. This arises from a parabolic rescaling and the central limit theorem, and it leads to a gradient flow with a (formal) inner-product structure, or equivalently, to a linear operator $M_\rho$ in~\eqref{def:M-GF}. Other types of randomness lead to non-quadratic gradient-flow structures, as we now describe.

\medskip

A close inspection of the arguments of Section~\ref{subsec:alt-GF-formulation} shows that they hinge on the inequality
\begin{equation*}
\partial_t \E(\rho) \geq -\frac12 \|\partial_t \rho\|_{M^{-1}_\rho}^2 - \frac12 \Bigl\|-\frac{\delta \E}{\delta\rho}\Bigr\|^2_{M_\rho},
\end{equation*}
together with the observation that equality holds if and only if $\partial_t \rho = - M_\rho\delta\E/\delta\rho$. This can be generalised by introducing a Legendre pair of convex functions $\psi_\rho$ and $\psi_\rho^*$, where the subscript~$\rho$ serves to indicate that they may depend on $\rho$, in the same way as the operator~$M_\rho$ does; in this context, $\psi_\rho$, is often called dissipation potential. In terms of this pair we then derive that
\begin{equation*}
\partial_t \E(\rho(t)) = \int \frac{\delta\E}{\delta\rho}\partial_t \rho 
\geq -\psi_\rho^*(\partial_t\rho) - \psi_\rho\Bigl(\frac{\delta\E}{\delta\rho}\Bigr),
\end{equation*}
and equality holds if and only if 
\begin{equation}
\label{def:genGF}
\partial_t \rho \in \partial \psi_\rho\Bigl(- \frac{\delta\E}{\delta\rho}\Bigr).
\end{equation}
%Here $\partial\psi_\rho^*$ is the subdifferential of $\psi_\rho^*$. 
The case of the $M$-gradient flow~\eqref{def:genGF} corresponds to 
\begin{equation*}
\psi_\rho^*(\xi) := \frac12 \|\xi\|^2_{M_\rho}
\qquad\text{and}\qquad
\psi_\rho(s) := \frac12 \|s\|_{M^{-1}_\rho}^2 .
\end{equation*}

The obvious generalisation of~\eqref{equiv:connection} then is
\begin{equation}
\label{equiv:connection-psi}
\text{$\rho$ is a solution of the $(\E,\psi)$-gradient-flow~\eqref{def:genGF}} \quad\Longleftrightarrow \quad J_\psi(\rho)=0, 
\end{equation}
where $J_\psi$ is given by
\begin{equation}
\label{def:Jpsi}
J_\psi(\rho) := \E(\rho(T)) - \E(\rho(0)) + \int_0^T\left[\psi_\rho^*(\partial_t \rho) + \psi_\rho\Bigl(- \frac{\delta\E}{\delta\rho}\Bigr)\right]\, dt.
\end{equation}

\subsection{Birth-death processes}

A simple example of a stochastic process with non-quadratic dissipation $\psi$ and a corresponding generalised gradient flow is a birth-death process, which is a continuous-time jump process on  $\Z$. The system may only jump to neighbours, from position $k$ with rate $a_k$ to $k+1$ and with rate $b_k$ to $k-1$. We construct a continuum limit by defining the new stochastic variable $U_n$ by rescaling time~$t$ and position $k(t)$ with $n$:
\begin{equation*}
U_n(t) := \frac{k(nt)}n.
\end{equation*}
A standard argument gives the large-deviation behaviour for $U_n$ in terms of the rate functional (see~\cite{Chan98} for a finite-lattice proof of the claims made below). 
If we choose the jump rates so that 
\begin{equation*}
a_k = \alpha e^{-\E'(k/n)}\qquad\text{and}\qquad
b_k = \alpha e^{+\E'(k/n)}
\end{equation*}
for $\alpha>0$ and some smooth function $\E\colon\R\to\R$, then the rate functional is 
\begin{equation*}
I(u) = \int_0^T L(u(t),u'(t))\, dt,
\end{equation*}
with
\begin{equation*}
L(u,v) = v\log \frac{v + \sqrt{v^2 + 4\alpha^2}}{2\alpha\exp(-\E'(u))}
- \sqrt{v^2 + 4\alpha^2} + \alpha e^{-\E'(u)} + \alpha e^{+\E'(u)}.
\end{equation*}
Writing
\begin{equation*}
\psi^*(v) = v\log \frac{v+\sqrt{v^2+4\alpha^2}}{2\alpha} - \sqrt{v^2+4\alpha^2},
\end{equation*}
it follows that $\psi(\xi) = \alpha(e^\xi+e^{-\xi})$, and $I$ can be written in the form~\eqref{def:Jpsi}. 

The corresponding generalised gradient flow in $\R$, given by~\eqref{def:genGF}, reads
\begin{equation*}
\dot u = 2\sinh (-\E'(u)).
\end{equation*}
Observe how this differs from the standard (quadratic-dissipation) gradient flow, which is $\dot u = -\E'(u)$; the non-quadratic dissipation preserves the sign of the velocity, but not its amplitude. Because of the preservation of sign, the energy $\E$ is monotonic along a solution:
\begin{equation*}
\frac{d}{dt} \E(u(t)) = \E'(u(t))\dot u(t) = 2\E'(u)\sinh (-\E'(u))\leq 0.
\end{equation*}

This example shows how the connection between large-deviation principles and (generalised) gradient flows extends to the case of non-quadratic dissipations. Note that here the large deviations refer to a single process and henceforth are not due to an averaging process as in the empirical measure case.

\subsection{Spin-flip processes}

For $ n\in\mathbb{N} $, let $ \mathbb{T}_n $ be the one-dimensional $ n$-torus $ (\mathbb{Z}/n\mathbb{Z}) $. An Ising spin at  sites of $ \mathbb{T}_n $ takes values in $ \{-1,+1\} $ and is subject to a rate-$1$ independent spin-flip dynamics.  We consider the trajectory of the magnetisation, i.e., $ t\mapsto m_n(t)=\frac{1}{n}\sum_{i\in\mathbb{T}_n}\sigma_i(t)$, where $ \sigma_i(t)$ is the spin at site $ i\in\mathbb{T}_n $ at time $t$. The generator for the process $ (m_n(t))_{t\ge 0} $ is given by
\[
(A_n f)(m)=\frac{(1+m)}{2}n[f(m-2n^{-1})-f(m)]+\frac{(1-m)}{2}n[f(m+2n^{-1})-f(m)]
\]
for $m\in\{-1,-1+2n^{-1},\ldots,1\} $. The trajectory of the magnetisation satisfies a large deviation principle, i.e., for every trajectory $ \gamma=(\gamma_t)_{t\in[0,T]} $,
\[
\Prob\big{(}(m_n(t))_{t\in[0,T]}\approx(\gamma_t)_{t\in[0,T]}\bigr{)}\approx\exp\Big{[}-n\int_0^TL(\gamma_t,\dot{\gamma}_t) \,dt\Big{]},
\]
where the Lagrangian $ L $ can be computed following the scheme of Feng and Kurtz~\cite[Example~1.5.]{FengKurtz06}.  
We obtain
\[
L(m,q)=\frac{q}{2}\log\bigg{(}\frac{q+\sqrt{q^2+4(1-m^2)}}{2(1-m)}\bigg{)}-\frac{1}{2}\sqrt{q^2+4(1-m^2)} +1. 
\]
This can similarly be written as $\psi^*(q) + \psi(-\E^\prime(m)) + q\E^\prime(m) $, where
\[
\psi^*(q) = \frac q2 \log \frac{q+\sqrt{q^2 + 4(1-m^2)}}{2\sqrt{1-m^2}}  - \frac12 \sqrt{q^2 + 4(1-m^2)}
\]
and
\[
\psi(\xi) = \frac12 \sqrt{1-m^2}\; \bigl{(}\exp(2\xi) + \exp(-2\xi)\bigr{)};
\]
the involved energy is
\[
\E(m) = \frac14 (1+m)\log (1+m) + \frac14 (1-m)\log (1-m).
\]
Then the limiting equation~\eqref{def:genGF} can be written as $\dot m = -2m$. This is consistent with the optimal trajectory via the Euler-Lagrange equation, $ m(t)=m_0e^{-2t}$.

\section{Discussion}
\label{sec:discussion}

In the sections above we have described a number of  pairs of systems, each consisting of a stochastic process and its continuum limit. Each pair has the property that the large deviations of the stochastic process are closely linked to a gradient-flow structure of the limit equation. These links are time-dynamic versions of the equilibrium connection mentioned in the introduction.  We now describe how this provides us insight into the properties of the gradient-flow structures for each pair.

\subsection{Wasserstein gradient flows}
\label{sec:Wass-grad-flows}

We claim that the Wasserstein metric characterises the mobility of the empirical measure of a large number of Brownian particles.
Indeed, this claim can  be made meaningful in a number of different ways:
\begin{enumerate}
\item In discrete time, letting $\rho_n$ be the empirical measure of a system of Brownian particles,  we have 
\begin{equation*}
\Prob (\rho_n(h)\approx \rho^1|\rho_n(0)\approx \rho^0) \sim e^{-nh \cdot d(\rho^0,\rho^1)^2/4}
\qquad \text{as } n\to\infty,
\end{equation*}
which follows from~\eqref{def:K_h} and was proved independently in~\cite{Leonard07TR}.
\item In continuous time, for the whole path $\rho_n\colon [0,T]\to\mathcal{M}_1(\mathbb{R}^d)$ of empirical measures up to a fixed terminal time $T$, we have
\begin{equation*}
\Prob (\rho_n\approx \rho) \sim e^{-nI(\rho)}, 
\qquad \text{as } n\to\infty,
\end{equation*}
where $I$, defined in~\eqref{I:KipnisOlla}, measures the size of the deviation by the Wasserstein metric tensor $\|\cdot\|_\rho$.
\item When the particles also undergo a deterministic drift, the same statement holds with $I$ defined by~\eqref{I:FK}, where again the size of the deviation is measured by the norm $\|\cdot\|_\rho$.
\end{enumerate}
The origin of this role of the Wasserstein metric as the mobility of Brownian particles can be understood by considering the geometric relationship between $(\R^d)^n$ and the space of measures endowed with the Wasserstein distance. 
Consider the embedding
\begin{equation*}
e\colon(\R^d)^n \to (\M_1(\R^d),d), \qquad
(x_1,\dots,x_n) \mapsto \frac1n\sum_{i=1}^n \delta_{x_i}.
\end{equation*}
Note that $e$ is not one-to-one, since the numbering of the particles is lost: the particles have become indistinguishable. Indeed, one can identify the set of empirical measures of the form $n^{-1} \sum_i \delta_{x_i}$ with the space obtained by identifying all elements in $(\R^d)^n$ that are rearrangements of each other, i.e., the quotient space $(\R^d)^n/S_n$, where $S_n$ is the set of all permutations of $n$ elements. 

Now the Wasserstein metric on $\M_1$ makes the embedding of $(\R^d)^n/S_n$ in $\M_1(\R^d)$   \emph{isometric}. This follows from the simple property that 
\begin{equation}
\label{equiv:two_metrics}
d\Biggl{(}\frac1n \sum_{i=1}^n \delta_{x_i},\frac1n \sum_{j=1}^n \delta_{y_j}\Biggr{)}^2 
\;=\; \frac1n \inf_{\sigma\in S_n} \sum_{i=1}^n \bigl{|}x_i-y_{\sigma(i)}\bigr{|}^2.
\end{equation}
With this property the role of the Wasserstein distance can be fully explained. The Freidlin-Wentzell theory for Brownian particles~\cite{FreidlinWentzell98} shows how the mobility of a \emph{vector} $X=(X_1,\dots,X_n)$ of $n$ Brownian particles has a stochastic mobility given by the Euclidean norm $(x,y)\mapsto \sum_i |x_i-y_i|^2$, in the sense that 
\begin{equation*}
\Prob (X(h) \approx y|X(0) \approx x) \sim \exp -\frac 1{4h}  \sum_{i=1}^n |x_i-y_i|^2
\qquad\text{for small }h.
\end{equation*}
The loss of information upon introducing indistinguishability, or equivalently upon transforming to empirical measures, implies by the contraction principle (e.g.,~\cite[Sec. 4.2.1]{DemboZeitouni98}) that the exponent $(1/4)\sum_i |x_i-y_i|^2$ becomes replaced by its minimum under rearrangement, 
\begin{equation*}
\frac 1{4h}\inf_{\sigma\in S_n} \sum_{i=1}^n |x_i-y_{\sigma(i)}|^2.
\end{equation*}
This expression is equal to $n/(4h)$ times~\eqref{equiv:two_metrics}. If we gloss over the  approximations in different limits ($h\to0$ and $n\to\infty$), this explains how the Wasserstein distance is the natural measure of the mobility of an \emph{empirical measure} of Brownian particles, through transformation of the original mobility of a single Brownian particle.

\subsection{Consequences for modelling}

Gradient flows can be thought of as overdamped systems, in the sense that any inertial effects are damped out quickly by the effects of viscous, frictional, or other damping forces, and can therefore be neglected. 
One way of modelling such overdamped systems is therefore by assuming an abstract gradient-flow structure from the start and making it concrete by postulating an energy $\E$ and a dissipation potential $\psi$. These choices should be motivated, and in the case of  Wasserstein and Wasserstein-like dissipations this motivation is non-trivial.

One area where this is particularly visible is in the modelling of lower-dimensional structures, such as threads and surfaces, moving through a viscous fluid. The biology of sub-cell structures knows many such examples, including microtubules and lipid bilayers. The assumption of overdampedness is reasonable in this viscosity-dominated situation, but the interplay of geometry and mechanics makes the direct formulation of evolution equations complicated and error-prone (see, e.g.,~\cite{ArroyoDeSimone09}). In this context, the construction of evolution equations through the postulation of energy and dissipation is often simpler and allows for clearer separation of the various assumptions. However, it remains necessary to motivate the choices made for the energy and the dissipation.

To take the Wasserstein metric as an example, its interpretation as the measure of mobility of empirical measures of Brownian particles provides such a motivation, and because of the connection to the Brownian mobility of the particles it also allows for generalisation to other situations. 

But similar arguments apply to other dissipations, coupled to other underlying stochastic processes. For instance,  the symmetric simple exclusion process leads to $\rho(1-\rho)$ mobility, implying that if such an exclusion process is one's idea of the underlying system, then the $\rho(1-\rho)$-dissipation is the natural choice. 

One might go even further. The diffusion equation~\eqref{eq:diffusion} is known to be a gradient flow in many different ways; in addition to the two mentioned above, also as the $L^2$-gradient flow of the Dirichlet integral $\int|\nabla \rho|^2$, for instance,  as the $H^{-1}$-gradient flow of the $L^2$-norm, and even as the $H^{s-1}$-gradient flow of the $H^s$-seminorm for each $s\in\R$. For the two structures that we have discussed, the different underlying stochastic processes provide clear reasons for the differing dissipations and energies. Here we formulate the 
\begin{conjecture}
\textbf{Each} gradient-flow structure can be connected to an appropriate stochastic process via a large-deviation principle.
\end{conjecture}
To the extent that this conjecture turns out to be true, it provides an explanation for the occurrence of multiple gradient-flow formulations of the same differential equation.

\subsection{Geometry and reversibility}

\def\diffc{A}
There are interesting connections between the geometry of the Brownian noise, the reversibility of the stochastic process, and the question whether the resulting evolution equation is a gradient flow or not.

This becomes apparent when we modify the system of Section~\ref{subsec:interactions} by introducing a diffusion matrix $\diffc\in \R^{d\times d}$ and replacing the scalar $\sigma$ by a mobility matrix $\sigma\in \R^{d\times d}$, thus obtaining
\begin{equation}
\label{eq:SDE-disc}
d X_i(t) = -\diffc\nabla \Psi(X_i(t))\,dt - \frac1n \sum_{j=1}^n \diffc\nabla\Phi(X_i(t)-X_j(t))\, dt + \sqrt 2\, \sigma\, dW_i(t).
\end{equation}

The large-deviation rate functional  of the system is similarly given by 
\begin{equation}
\label{I:FK-disc}
I(\rho) := \frac12 \int_0^T \Bigl\|\partial_t\rho - \div \sigma\sigma^T\nabla \rho 
- \div \rho \diffc\nabla \bigl[\Psi + \rho\ast \Phi\bigr]\Bigr\|^2_{D(\rho),*}\, dt,
\end{equation}
where the norm $\|\cdot\|_{D(\rho),*}$ is induced by~\eqref{def:gD} with $D(\rho) = \rho\sigma\sigma^T$. 
The formula~\eqref{I:FK-disc} implies that the hydrodynamic limit of this system is the minimiser of $I$, satisfying
\begin{equation}
\label{eq:interaction-disc}
\partial_t\rho =  \div \sigma\sigma^T\nabla \rho 
+ \div \rho \diffc \nabla \bigl[\Psi + \rho\ast \Phi\bigr].
\end{equation}
%is a Wasserstein gradient flow of the functional
%\begin{equation*}
%E(\rho) := \int_{\R^d}\left[ \frac12 \rho\log\rho + \rho\Psi + \frac12 \rho(\rho\ast \Phi)\right].
%\end{equation*}
%and we now sketch how this limit can be derived with the new approach.

With this additional parameter freedom, it is not always possible to write~\eqref{I:FK-disc} in the form~\eqref{def:JM}. This depends on whether the cross term in~\eqref{I:FK-disc} is an exact differential, i.e., whether there exists a functional $\E$ such that
\begin{equation*}
\Bigl(\partial_t\rho,-\div \sigma\sigma^T\nabla \rho 
- \div \rho \diffc \nabla \bigl[\Psi + \rho\ast \Phi\bigr]\Bigr)_{D(\rho),*} = \partial_t \E(\rho).
\end{equation*}
This is the case if and only if $\sigma\sigma^T$ is a positive multiple of $\diffc$, a condition that is familiar from the fluctuation-dissipation theorem. In that case, and writing $\sigma\sigma^T= kT\diffc$ for some `temperature' $T>0$ and the Boltzmann constant $k$, 
\begin{equation*}
-\div \sigma\sigma^T\nabla \rho 
- \div \rho \diffc \nabla \bigl[\Psi + \rho\ast \Phi\bigr]
= M_\rho \frac{\delta\F}{\delta\rho},
\end{equation*}
where $M_\rho\xi$ is defined as $-\div D(\rho)\nabla\xi$ and the free energy $\F$ is a modification of~\eqref{def:free-energy-F},
\begin{equation*}
\F(\rho) := \Ent(\rho) + \frac1{kT}\int_{\R^d}\left[\rho\Psi + \frac12 \rho(\rho\ast \Phi)\right].
\end{equation*}
Then the rate functional $I$ can be written in the form~\eqref{def:J} as 
\begin{equation*}
\label{eq:FK-FDT-disc}
I(\rho) 
= \F(\rho(T))- \F(\rho(0)) + \frac12 \int_0^T 
\left[ \|\partial_t \rho\|^2_{D(\rho),*} 
 + \Bigl\|\frac{\delta \F}{\delta\rho}\Bigr\|^2_{D(\rho)}
\right] \, dt
\end{equation*}
and the evolution equation~\eqref{eq:interaction-disc} is the (modified, $D$-) Wasserstein gradient flow of~$\F$.

Our freedom to choose $\diffc$ and $\sigma$ separately gives us the insight that for this system the following four statements are equivalent:
\begin{enumerate}
\item $\sigma\sigma^T= kT\diffc$ for some $T>0$;
\item The evolution~\eqref{eq:interaction-disc} is a $D(\rho)$-Wasserstein gradient flow of $\F$;
\item The rate functional $I$ can be written in the form~\eqref{def:J};
\item For any  finite number $n$ of particles, the system~\eqref{eq:SDE-disc} is reversible.
\end{enumerate}
We expect that such an equivalence property, including the reversibility of the microscopic system, might hold  more generally.

\subsection{Diffusion with decay}

Yet another generalisation concerns systems with decay, which is implemented as a jump process.
In~\cite{PeletierRenger11TR}, Peletier and Renger have derived a similar connection for the case of diffusing particles that are convected and may also decay, given by the equation (in one space dimension)
\begin{equation}
\label{eq: diffusion drift decay}
	\partial_t \rho = \partial_{xx} \rho + \partial_x \!\left( \rho \, \partial_x\!\Psi \right) - \lambda \rho, \\
\end{equation}
with $\Psi\in C_b^2(\R)$ and $\lambda\geq 0$. 

In~\cite{PeletierRenger11TR}, the particles perform a Brownian motion in the spatial dimension, augmented by a deterministic drift given by $-\partial_x \Psi$. This part of the process gives rise to the two terms $\partial_{xx} \rho + \partial_x \!\left( \rho \, \partial_x\!\Psi \right)$. In addition, the particles  change their state from `normal' to `decayed', after an exponentially distributed time; this part gives rise to the term $-\lambda \rho$. The opposite transition is not allowed: decay is irreversible. 

An analysis similar to Section~\ref{subsec:just-diff} then connects the large-deviation rate functional for this stochastic particle system to a corresponding minimisation problem describing the time-discrete evolution, i.e., the equivalent of~\eqref{approx:JKO}. In this case the time-discrete minimisation problem is
\def\ND{_{N\!D}}
\def\Dec{_{Dec}}
\begin{multline}
\label{min:intro-rhok-DD}
\rho^k \in\argmin_{\rho} \ \inf_{\rho\ND : |\rho+\rho\ND|=|\rho^{k-1}|} 
-\frac12 \F(\rho+\rho\ND) -\frac12 \F(\rho^{k-1}) 
  + \frac{1}{4h} d(\rho +\rho\ND,\rho^{k-1})^2\\ 
+ \F(\rho) + \F(\rho\ND) -|\rho|\log e^{-\lambda h} - |\rho\ND|\log (1-e^{-\lambda h}),
\end{multline}
where $|\rho|:=\int \rho$ and the free energy $\F$ is defined as
\begin{equation*}
\F(\rho ) = \Ent(\rho) + \int \Psi\, d\rho.
\end{equation*}
%To interpret~\eqref{min:intro-rhok-DD}, one should realise that the infimum over the measure $\rho\ND$ represents a choice: in each time step, the system designates a portion $\rho\ND\geq 0$ for decay (the index $N\!D$ stands for `Normal to Decayed'), while the other part $\rho\geq 0$ remains `normal'. 

In~\cite{PeletierRenger11TR}, the authors explain how the structure of~\eqref{min:intro-rhok-DD} can be understood: if we define 
\def\init#1{\overline{#1}}
\begin{align*}
K^h_\Psi(\init\rho; \rho^{k-1}) \;&:=\; \frac12\F(\init\rho) - \frac12\F(\rho^{k-1}) + \frac1{4h} d(\init \rho,\rho^{k-1})^2 ,  \\
K^h\Dec(\rho;\init \rho) \;&:=\;  \F(\rho) -\F(\init \rho) + \F(\init\rho\!-\!\rho) \;-\; |\rho|\log e^{-\lambda h} \;+\; |\init\rho-\rho|\log (1-e^{-\lambda h}) , 
\end{align*}
then the terms inside the infimum in~\eqref{min:intro-rhok-DD} can be written as $K^h(\rho+\rho\ND;\rho^{k-1})+ K^h\Dec(\rho;\rho+\rho\ND) $.  
In this decomposition, the first term describes diffusion and convection by $\Psi$ of the joint measure $\rho+\rho\ND$ starting from the previous state $\rho^{k-1}$,  similar to~\eqref{def:K_h} and~\eqref{approx:JKO}. The second term describes the decay process, in which the joint diffused-and-convected measure $\rho+\rho\ND$ is split into a part $\rho$ that remains `normal' and the remainder $\rho\ND$ that becomes decayed. 

While the structure of~\eqref{min:intro-rhok-DD} is not the same as~\eqref{approx:JKO}, and~\eqref{min:intro-rhok-DD} does not represent a time discretisation of a gradient flow, both are minimisation problems that define the next step in the iteration, and in both cases one can identify a driving force (the free energy $\F$, in the case of~\eqref{min:intro-rhok-DD}) and a mechanism that acts as a brake. In $K^h_\Psi$ the `brake' is the Wasserstein metric $d(\init\rho,\rho^{k-1})^2/4h$, and in $K^h\Dec$ it is the two terms $-|\rho|\log e^{-\lambda h} \;+\; |\init\rho-\rho|\log (1-e^{-\lambda h})$. In both cases these terms restrict the movement of $\init\rho$ respectively $\rho$, and this restriction becomes more and more severe as $h\to0$.

\subsection{General remarks on interacting particle systems}

Section~\ref{sec:Other-gradient-flows} explained how, once a large deviation principle for the interacting particle system with rate functional $I(\rho)$ is established, different Wasserstein-type metrics occur in a natural way. Such large deviation results are stronger than results on limit equations. Indeed, a part of the standard proof of a large deviation result involves modifying the process by adding a forcing such that a given path which does not solve the original limit equation solves the limit equation of the modified process. So the question arises whether the point of view advocated in this paper has the potential of deriving limit equations without using large deviation results which contain limit results derived in the classical way. This open question is of particular importance because limit points of the implicit time discretisation provide a weak notion of solution of the limit equation in cases where distributional solutions are not appropriate, e.g., for problems with a sharp interface like the mean curvature flow. In situations such as~\eqref{eq:SDE}, where a particle interacts with the average of many others, the distribution of a family of initially independent particles stays close to a product measure (propagation of chaos), so a modification of the techniques for independent particles seems promising.

\section{Conclusion}

The examples of this paper illustrate how the two concepts of large-deviation principles for stochastic particle systems and gradient flows are closely entwined. Further examples are currently under study, such as Brownian particles with inertia, which lead to the Kramers' equation, and rate-independent systems such as  friction and fracture. 
We expect that many more examples of this kind will be uncovered.  

\appendix

\section{Free energy and the Boltzmann distribution}
\label{subsec:free-energy}

In this appendix we show how the \emph{free energy}
\begin{equation}
\label{def:app:free-energy}
\mathcal F(\rho) := H(\rho|\mu) + \frac1{k\theta} E(\rho)
\end{equation}
arises from the coupling of a system of particles with a \emph{heat bath}. Here $\theta>0$ (in Joules) is the temperature of the heat bath, and the \emph{Boltzmann constant} $k$ has the value $1.4 \cdot 10^{-23}\un{J/K}$. The  measure $\mu\in \P(\X)$ is the probability distribution of the particles in a state space $\X$, and~$E$ is the average energy of the particles:
%arises from the coupling of a system of particles with a \emph{heat bath}. Here $\theta>0$  is the temperature of the heat bath, and $k$ is the \emph{Boltzmann constant}. The  measure $\mu\in \P(\X)$ is the probability distribution of the particles in a state space $\X$, and~$E$ is the average energy of the particles:
\begin{equation*}
E(\rho) = \int_{\X } e(x)\,\rho(dx),
\end{equation*}
where $e\colon\X\to\R$ is a fixed function that we call the \emph{energy} of a state $x\in \X$.
We now construct an explicit system in which $\F$ arises as the large-deviation rate functional. This will allow us to interpret all these concepts in the context of large deviations. 

\medskip
We start by  choosing a system $S$ and its connection to a heat bath called $S_B$. Both are probabilistic systems of particles; $S$ consists of $n$ independent particles $X_i\in\X$, with probability law $\mu\in\P(\X)$; similarly $S_B$ consists of $m$ independent particles $Y_j\in\Y$, with law $\nu\in\P(\Y)$. The total state space of the system is therefore $\X^n\times\Y^m$.

The \emph{coupling} between these systems is done  via an \emph{energy constraint}.  We assume that there are energy functions $e\colon\X\to\R$ and $e_B\colon\Y\to\R$, and we will constrain the joint system to be in a state of fixed total energy, i.e., we will only allow states in $\X^n\times\Y^m$ that satisfy 
\begin{equation}
\label{eq:energy-constraint}
\sum_{i=1}^n e(X_i) + \sum_{j=1}^m e_B(Y_j) = \text{constant}.
\end{equation}
The physical interpretation of this is that energy (in the form of heat) may flow freely from one system to the other, but no other form of interaction is allowed. 

Similar to the example in the Introduction, we describe the total states of systems $S$ and $S_B$ by empirical measures $\rho_n = \frac1n \sum_i \delta_{X_i}$ and $\zeta_m = \frac1m\sum_j \delta_{Y_j}$. We define the average energies $E(\rho_n) := \frac1n\sum_i e(X_i) = \int_{\X}e\, d\rho_n$ and $E_B(\zeta_m) := \int_{\Y} e_B\,d\zeta_m$, so that the energy constraint~\eqref{eq:energy-constraint} reads $nE(\rho_n) + mE_B(\zeta_m) = \text{constant}$.  

By Sanov's theorem each of the systems \emph{separately} satisfies a large-deviation principle with rate functions $I(\rho) = \H(\rho|\mu)$ and $I_B(\zeta) = \H(\zeta|\nu)$. However, instead of using the explicit formula for $I_B$, we are going to assume that $I_B$ can be written as a function of the energy $E_B$ of the heat bath alone, i.e., $I_B(\zeta) = \tilde I_B(E_B(\zeta))$. For the coupled system we derive a joint large-deviation principle by choosing that (a) $m = nN$ for some large $N>0$, and (b) the constant in~\eqref{eq:energy-constraint} scales as $nN$, i.e., 
\begin{equation*}
nE(\rho_n) + nNE_B(\zeta_{nN}) = nN\overline E \qquad\text{for some }\overline E.
\end{equation*}
Formally, the joint system then satisfies  a large-deviation principle
\begin{equation*}
\Prob\Big((\rho_n,\zeta_{nN}) \approx (\rho,\zeta) \;\Big|\; E(\rho_n) + NE_B(\zeta_{nN}) = \overline E\Bigr) \sim \exp \bigl(- nJ(\rho,\zeta)),
\end{equation*}
with rate functional
\begin{equation*}
J(\rho,\zeta) := 
\begin{cases}
\H(\rho|\mu) +N \tilde I_B(E_B(\zeta))+\text{constant}& \text{if } E(\rho) + NE_B(\zeta) = N\overline E,\\
+\infty&\text{otherwise.}
\end{cases}
\end{equation*}
Here the constant is chosen to ensure that $\inf J=0$.

The functional $J$ can be reduced to a functional of $\rho$ alone,
\begin{equation*}
J(\rho) = \H(\rho|\mu) + N\tilde I_B\left(\overline E - \frac{ E(\rho)}N\right)+\text{constant}.
\end{equation*}
In the limit of large $N$, one might approximate 
\begin{equation*}
N\tilde I_B\left(\overline E - \frac{E(\rho)}N\right) \approx N\tilde I_B(\overline E) - \tilde I_B'(\overline E)E(\rho).
\end{equation*}
The first term above is absorbed in the constant, and we find
\begin{equation*}
J(\rho) \approx \H(\rho|\mu) -E(\rho)\tilde I_B'(\overline E) +\text{constant}.
\end{equation*}
We expect that $I_B'$ is negative, since larger energies typically lead to higher probabilities and therefore smaller values of $I_B$. 
Now we simply define $k\theta := -1/\tilde I_B'(\overline E)$, and we find 
\begin{equation*}
J(\rho) \approx \H(\rho|\mu)+\frac1{k\theta} E(\rho) +\text{constant}.
\end{equation*}
This is the same expression as~\eqref{def:app:free-energy}.
Note that the right-hand side can be written as $\H(\rho|\tilde \mu)$, where $\tilde \mu$ is the tilted distribution
\begin{equation*}
\tilde \mu(A) = \frac{\ds\int_A e^{-e(x)/k\theta}\, \mu (dx)}{\ds\int_{\X} e^{-e(x)/k\theta}\, \mu (dx)}.
\end{equation*}

This derivation  shows that the effect of the heat bath is to \emph{tilt} the system $S$: a state $\rho$ of~$S$ with larger energy $E(\rho)$  implies a smaller energy $E_B$ of $S_B$, which in turn reduces the probability of~$\rho$. This is reflected in the approximation $\tilde I'_B(\overline E) E(\rho)$ of $I_B(\zeta)$. The role of temperature $\theta$ is that of an exchange rate, since it characterises the change in probability (as measured by the rate function $I_B$) per unit of energy. When $\theta$ is large, the exchange rate is low, and then larger energies incur only a small probabilistic penalty. When temperature is low, then higher energies are very expensive, and therefore more rare.
From this point of view, the Boltzmann constant $k$ is simply the conversion factor that converts our Kelvin temperature scale for $\theta$ into the appropriate `exchange rate' scale.

In  thermodynamics one often encounters the identity (or definition) $\theta = dS/dE$. This is formally the same as our definition of $k\theta$ as $-dI_B/dE$, if one interprets $I_B$ as an entropy and adopts the convention to multiply the non-dimensional quantity $I_B$ with $-k$. 
%I find this latter way of thinking less intuitive than the one explained above, however. 

\medskip

%Yet another insight that this calculation gives is the following. I was often puzzled by the fact that in defining a free energy (e.g., $E-\theta S$, or $k\theta \Ent +E$), one adds two rather different objects: the energy of a system seems to be a completely different type of object than the entropy. This derivation of Boltzmann statistics shows that  it's not exactly energy and entropy that one's adding; it really is more like adding two entropies ($\H$ and $I_B$, in the notation above). The fact that we write the second entropy as a constant times energy follows from the coupling and the approximation allowed by the assumption of a large heat bath.
%the energy is acting as a tilting of the original behaviour (described by the entropy of the original system), and therefore it's the original entropy
%

%One might ask how this interpretation of temperature matches with the experience that temperature determines thermal equilibrium (two bodies are in thermal equilibrium if and only if their temperatures are equal) and heat flux (in the direction of decreasing temperature). To some extent this can be explained by a similar argument. We think of two bodies being in thermal equilibrium if `nothing happens' when we allow them to exchange heat, as follows. We consider two systems with energies $E_1$ and $E_2$ in equilibrium (whatever that may mean) and we couple them by fixing the sum $E_1+E_2=\overline E$, while allowing for variations in $E_1$ and $E_2$---very much as above. The systems are then said to be in thermal equilibrium if the probability distributions before and after

\bigskip

\noindent
\textbf{Acknowledgement. } The authors wish to thank Dejan Slep\v cev and Rob Jack for various interesting discussions. The research of M.~A.~Peletier and J.~Zimmer has received funding from the Initial Training Network
``FIRST'' of the Seventh Framework Programme of the European Community
(grant agreement number 238702). The research of S. Adams was supported by EPSRC grant number EP/I003746/1.

\bibliographystyle{plain}
\bibliography{ref-gf}
\end{document}